\documentclass[12pt]{amsart}
\usepackage{latexsym, amssymb, amsmath, amscd}
 \usepackage{eucal}

    \newtheorem{rema}{Remark}[section]
    \newtheorem{propo}[rema]{Proposition}

   \newtheorem{theo}[rema]{Theorem}
   \newtheorem{def-theo}[rema]{Definition-Theorem}
 \newtheorem{conj}[rema]{Conjecture}
   \newtheorem{defi}[rema]{Definition}
    \newtheorem{lemma}[rema]{Lemma}
    \newtheorem{corol}[rema]{Corollary}
     \newtheorem{exam}[rema]{Example}

    \newtheorem{prob}[rema]{Problem}

	\newcommand{\nno}{\nonumber}

 \newcommand{\pf}{{\it Proof:}\hspace{2ex}}
 
 \newcommand{\epfv}{\hspace{1em}$\Box$\vspace{1em}}

\newcommand{\bC}{{\mathbb C}}
\newcommand{\bZ}{{\mathbb Z}}

\newcommand{\bN}{{\mathbb N}}

\newcommand{\cL}{{\mathcal L}}

\newcommand{\cA}{{\mathcal A}}

\newcommand{\Ker}{\mbox{\rm Ker\,}}

 \newcommand{\rg}{ {R[G]} }
 \newcommand{\kg}{ {K[G]} } 
\newcommand{\tr}{ {\rm Const\hspace{0.3mm}} } 
\newcommand{\vg}{ {V_G} } 
\newcommand{\kz}{K[z]}
\newcommand{\rz}{R[z]} 
\newcommand{\lin}{1\le i\le n}

\title[Analogue of the Duistermaat-van der Kallen Theorem]
{Analogue of the Duistermaat-van der Kallen Theorem 
for Group Algebras}

\author{Wenhua Zhao and Roel Willems}

\address{W. Zhao, 
 Illinois State University,  Normal, IL  
61790-4520, USA.
{\em Email}: wzhao@ilstu.edu }

\address{R. Willems, 
Radboud University Nijmegen,  
Postbus 9010, 6500 GL Nijmegen, 
The Netherlands.   
{\em Email}: r.willems@math.ru.nl}

  \date{\today} 

\begin{document}

\begin{abstract} 
Let $G$ be a group, $R$ an integral domain, and $\vg$ 
the $R$-subspace of the group algebra $\rg$ consisting of 
all the elements of $\rg$ whose coefficient of the identity element 
$1_G$ of $G$ is equal to zero. Motivated by 
the Mathieu conjecture \cite{Ma}, the Duistermaat-van der Kallen theorem  \cite{DK}, and also by recent studies on the notion of Mathieu subspaces introduced 
in \cite{GIC} and \cite{MS}, we show that for finite groups $G$, 
$\vg$ under certain conditions also forms a Mathieu subspace of the group algebra $\rg$. We also show that for the free abelian 
groups $G=\bZ^n$ $(n\ge 1)$ and any integral domain $R$ of positive characteristic, $\vg$ fails to be a Mathieu subspace of $\rg$, which is equivalent to saying that the Duistermaat-van der Kallen 
theorem  \cite{DK} cannot be generalized to any 
field or integral domain of positive characteristic.    
\end{abstract} 

\keywords{Duistermaat-van der Kallen Theorem, 
Mathieu subspaces, groups algebras} 
   
\subjclass[2000]{16S34, 16N40, 16D99}

\thanks{The first author has been partially 
supported by NSA Grant H98230-10-1-0168 and the second author 
is funded by Phd-grant of council for the physical sciences, 
Netherlands Organization for scientific 
research (NWO)}

 \bibliographystyle{alpha}
    \maketitle


\renewcommand{\theequation}{\thesection.\arabic{equation}}
\renewcommand{\therema}{\thesection.\arabic{rema}}
\setcounter{equation}{0}
\setcounter{rema}{0}
\setcounter{section}{0}

\section{\bf Introduction}
Let's first recall the following notion 
introduced recently by the first author  
in \cite{GIC} and \cite{MS}, which can be viewed as 
a natural generalization of the notion of ideals. 
 
\begin{defi}\label{Def-MS}
Let $R$ be a commutative ring and $\cA$ 
an associative $R$-algebra. A $R$-submodule or 
$R$-subspace $M$ of $\cA$ is said to be a {\it left} 
$($resp., {\it right}; {\it two-sided}$)$ {\it  Mathieu subspace} 
of $\cA$ if for any $a, b, c\in \cA$ with 
$a^m\in M$ for all $m\ge 1$, we have 
$ b a^m \in M$ $($resp., $a^m b \in M$; $ba^mc\in M$$)$ 
when $m\gg 0$, 
i.e., there exists $N\ge 1$ such that $b a^m  \in M$ 
$($resp., $a^m b \in M$; 
$ba^mc\in M$$)$ for all $m\ge N$.
\end{defi}   

Two-sided Mathieu subspaces will also simply be 
called Mathieu subspaces. A $R$-subspace $M$ of $\cA$ is said to be a {\it pre-two-sided} Mathieu subspace of $\cA$ if it is both left and right Mathieu subspace of $\cA$. Note that the {\it pre-two-sided} Mathieu subspaces were previously called  {\it two-sided} Mathieu subspace 
or {\it Mathieu subspaces} in \cite{GIC}. 

The introduction of the notion of Mathieu subspaces in \cite{GIC} 
and \cite{MS} was mainly motivated by the studies of 
{\it the Jacobian conjecture} \cite{K} 
(see also \cite{BCW} and \cite{E}), 
{\it the Mathieu conjecture}  \cite{Ma}, 
{\it the vanishing conjecture} \cite{HNP}, \cite{GVC}, \cite{AGVC}, 
\cite{EWiZ} and more recently, {\it the image conjecture} \cite{IC} 
as well as many other related open problems. 
For some recent developments on Mathieu subspaces, 
see \cite{MS}, \cite{FPYZ}, \cite{EWZ1}, \cite{EWZ2}, 
\cite{EZ} and \cite{GMS}.    
For a recent survey on the {\it the image conjecture} 
and it's connections with some other problems, 
see \cite{E2}.

The notion was named after Olivier Mathieu in \cite{GIC} due to his conjecture mentioned above, which now in terms of the new notion 
can be re-stated as follows.  

\begin{conj} \label{MC} $(${\bf The Mathieu Conjecture}$)$ 
Let $G$ be a compact connected real Lie group with the 
Haar measure $\sigma$. Let $\cA$ be the algebra of 
complex-valued $G$-finite functions on $G$, and 
$M$ the subspace of $\cA$ consisting of $f\in \cA$  
such that $\int_G  f \,  d\sigma=0$. Then 
$M$ is a Mathieu subspace of $\cA$. 
\end{conj}

J. Duistermaat and W. van der Kallen \cite{DK} proved 
{\it the Mathieu conjecture} for the case of tori, 
which now can be re-stated as follows.  

\begin{theo}\label{ThmDK} $(${\bf Duistermaat and 
van der Kallen}$)$  Let $z=(z_1, z_2,$ $..., z_n)$ 
be $n$ commutative free variables and 
$V$ the subspace of the Laurent polynomial algebra 
$\bC[z^{-1}, z]$ consisting of the Laurent polynomials 
with no constant term. Then $V$ is a Mathieu 
subspace of $\bC[z^{-1}, z]$. 
\end{theo}

Note that despite its innocent looking, the proof of 
the theorem above is surprisingly difficult. 
The proof in \cite{DK} uses some heavy machineries  
such as toric varieties, resolutions of singularities, 
etc. 

To discuss the main motivations and results 
of this paper, we start with the following observation 
on the Duistermaat-van der Kallen Theorem above. 

Let $G$ be the free abelian group $\bZ^n$ $(n\ge 1)$. 
Then the Laurent polynomial algebra $\bC[z^{-1}, z]$ 
can be identified in the obvious way 
with the group algebra $\bC[G]$.
Under this identification, the subspace 
$V\subset \bC[z^{-1}, z]$ in the theorem 
corresponds to the subspace $V_G$ 
of the group algebra $\bC[G]$ consisting of 
the elements of $\bC[G]$ 
whose ``{\it constant term}" 
(i.e., the coefficient of the identity element 
$1_G$ of $G$) is equal to zero. 
So, we are naturally led to the following 
(open) problem. 

\begin{prob}\label{MainProb} 
Let $R$ be a commutative ring and 
$G$ a group. Let $\vg$ be the $R$-subspace 
of the elements of the group algebra $R[G]$ 
with no ``constant term", i.e., 
the coefficient of the identity element 
$1_G$ of $G$ is equal to zero. 
Then under what conditions on $R$ and $G$, 
$\vg$ forms a Mathieu subspace of the 
group algebra $R[G]$? 
\end{prob}

The problem above not only provides 
a different point of view to get 
further understanding on the 
remarkable Duistermaat-van der 
Kallen Theorem, but also gives a family of 
candidates for Mathieu subspaces, which may provide 
some new understandings on the still very mysterious 
notion of Mathieu subspaces. This makes 
the problem itself very interesting  
and worthy to investigate.

One of the main results of this paper is 
that for any finite group $G$ and an  
integral domain $R$ of characteristic 
$p=0$ or $p>|G|$ (the {\it order} of $G$), 
the $R$-subspace $\vg$ does form a Mathieu 
subspace of $\rg$ (see Theorem \ref{GeneralCase}), 
i.e., Problem \ref{MainProb} in this case 
can be solved completely.  
 
However, for the case that $0<char.\,R=p \le |G|$, 
the situation becomes much more subtle.  
For example, the magic condition 
$p\nmid |G|$ for the group algebras 
of finite groups $G$ (e.g., see \cite{P}) 
does not resolve the difficulty completely, i.e., 
under this condition $\vg$ still may or 
may not be a Mathieu subspaces of 
$\rg$ (e.g., see Theorem \ref{MainAbThm} 
and Example \ref{F3Z5}). 

In this paper, we first study Problem \ref{MainProb}  
for the group algebras of finite groups $G$ 
over integral domains $R$ of any characteristics. 
In particular, besides the main result mention above, 
for finite abelian groups we also give a complete solution 
of Problem \ref{MainProb} for the case that the base 
integral domain $R$ satisfies certain primitive root of unity 
conditions (see Theorems \ref{GeneralCase} 
and \ref{MainAbThm}), e.g., when $R$ 
is an algebraically 
closed field.     
 
We then show that for the group algebras  
of the free abelian groups $G=\bZ^n$ $(n\ge 1)$  
over any integral domain $R$ of positive 
characteristic, $\vg$ is not a Mathieu subspace 
of $\rg$, by showing that an example suggested 
by Arno van den Essen does provide 
a desired counter-example. 
Consequently, it follows that 
the Duistermaat-van der Kallen theorem, Theorem \ref{ThmDK},  
cannot be generalized to the Laurent polynomial algebra 
$R[z^{-1}, z]$ over any field or integral domain $R$ of 
positive characteristic.  


The arrangement of this paper is as follows.


In Section \ref{S2}, we recall some general results on 
Mathieu subspaces obtained in \cite{GIC} and \cite{MS}, 
which will be needed later in this paper. 
In Section \ref{S3}, we prove some 
results on Problem \ref{MainProb} for the group algebras 
of finite groups $G$ over arbitrary commutative rings 
or integral domains. 
In particular, we show in Theorem \ref{GeneralCase} that 
when the base ring $R$ is an integral domain of 
characteristic $p=0$ or $p>|G|$, the subspace $\vg$ is 
always a Mathieu subspace of $\rg$.     

In Section \ref{S4}, we focus on the group algebras 
of finite abelian groups $G$ over integral domains $R$  
of characteristic $p>0$. The main results of this section 
is Theorem \ref{MainAbThm}, which combining with  
Theorem \ref{GeneralCase} provides a complete  
solution of Problem \ref{MainProb} for the group algebras 
of finite abelian groups $G$ over the integral domains $R$ 
which satisfies a primitive root of unity 
condition, e.g., when $R$ is an 
algebraically closed field.    

In Section \ref{S5}, we consider Problem \ref{MainProb} 
for the group algebras of the free abelian groups  
$\bZ^n$ $(n\ge 1)$ over an integral domain $R$ 
of characteristic $p>0$. We prove that $\vg$ in this case fails 
to be a Mathieu subspace of $R[\bZ^n]$ by showing 
that the example in Lemma \ref{ExampleArno}, which was suggested 
by Arno van den Essen to the authors, does provide  
a desired counter-example.

\renewcommand{\theequation}{\thesection.\arabic{equation}}
\renewcommand{\therema}{\thesection.\arabic{rema}}
\setcounter{equation}{0}
\setcounter{rema}{0}
 
\section{\bf Some Results on Mathieu Subspaces}\label{S2}

In this section, we recall some general facts on 
Mathieu subspaces which will be needed later in this paper. 
Although all the results below with certain modifications 
hold for all types of Mathieu subspaces 
(one-sided, pre-two-sided, etc.)  
We here only focus on the two-sided case, 
which by Corollary \ref{AllType} 
in the next section will be enough 
for our purpose. 

Throughout this paper, unless stated otherwise, 
$R$ and $K$ always stand respectively 
for a unital commutative ring and a field of 
any characteristic, and $\cA$ a unital algebra 
over $R$ or $K$.  

Following \cite{MS}, we define for any $R$-subspace 
$V$ of a $R$-algebra $\cA$ the {\it radical}, denoted by $\sqrt V$, 
to be the set of $a\in \cA$ such that $a^m\in V$ when $m\gg 0$.    

We start with the following equivalent formulation of Mathieu subspaces, which was given in Proposition $2.1$ 
in \cite{MS}.  

\begin{propo}\label{MS-Def2} 
Let $\cA$ be a $R$-algebra and $V$ a $R$-subspace 
of $\cA$. Then $V$ is a Mathieu subspace 
of $\cA$ iff for any $a\in \sqrt V$ and 
$b, c\in \cA$, we have $b a^m c \in V$ 
when $m\gg 0$.  
\end{propo}

The following characterization of the Mathieu subspaces 
with algebraic radicals was also proved in 
Theorem $4.2$ in \cite{MS}. 

\begin{theo}\label{CharByIdem}
Let $A$ be a $K$-algebra and 
$V$ a $K$-subspace of $\cA$ 
such that $\sqrt V$ is algebraic over $K$ 
$($i.e., every element of $\sqrt V$ is algebraic over $K$$)$.
Then $V$ is a Mathieu subspace of $\cA$ iff 
for any idempotent $e \in V$ $ ($i.e., $e^2=e$$)$, 
we have $(e) \subseteq V$, where $(e)$ denotes 
the ideal of $\cA$ generated by $e$. 
\end{theo}

The next proposition is easy to check directly 
(or see Proposition $2.7$ in \cite{MS}).    

\begin{propo}\label{quotient-I}
Let $I$ be an ideal of $\cA$ and 
$V$ a  $R$-subspace of $\cA$ 
such that $I\subseteq V$. 
Then $V$ is a Mathieu subspace of 
$\cA$ iff $V/I$ is a Mathieu subspace   
of the quotient algebra $\cA/I$.
\end{propo}

Finally, let's recall the following family of 
Mathieu subspaces of the polynomial algebra 
$\kz$ in $n$ variables $z\!:=(z_1, z_2, ..., z_n)$, 
which was given in Proposition $4.6$ in \cite{GIC}. 

\begin{propo}\label{FiniteCase2}
Let $n, d\ge1$ and $R$ an arbitrary integral domain.  
Let $S=\{v_1, v_2, ..., v_d \}\subset R^n$ 
$($with $d$ distinct elements$)$    
and $0 \ne c_i \in R$ $(1\le i\le d)$.  
Denote by $V$ the subspace of $f(z)\in \rz$ 
such that 
\begin{align}\label{FiniteCase2-e1}
\sum_{i=1}^d c_i f(v_i)=0.
\end{align} 
Then $V$ is a 
Mathieu subspace of $\rz$ iff 
for any non-empty subset 
$J\subset \{1, 2, ..., d\}$, we have
\footnote{Note that Eq.\,(\ref{FiniteCase2-e2}) in \cite{GIC} 
had been misprinted.}  
\begin{align}\label{FiniteCase2-e2}
\sum_{i\in J} c_i \ne 0. 
\end{align}
\end{propo}

Note that the proposition above was only proved 
in \cite{GIC} under the condition that $R$ is a field. 
But, it is easy to see that the same proof actually 
goes through equally well for all integral domains.

\renewcommand{\theequation}{\thesection.\arabic{equation}}
\renewcommand{\therema}{\thesection.\arabic{rema}}
\setcounter{equation}{0}
\setcounter{rema}{0}
 
\section{\bf Some General Results for the Case of Finite Groups}\label{S3}

Throughout the rest of this paper, unless stated otherwise, 
$G$ stands for a finite group, $R$ a commutative ring, and 
$K$ a field of any characteristic.  We denote by $R[G]$ 
and $\kg$ the group algebra of $G$ over $R$ and $K$, 
respectively.  Furthermore, we also fix the following 
terminologies and notations. 
 \begin{enumerate}
\item[$i)$] We denote by $1$ or $1_G$ 
the identity element of the group $G$  
and also the identity element of 
the group algebra $\rg$. 

\item[$ii)$]  For any $u\in \rg$, we denote by $\tr(u)$ 
the coefficient of $1_G$ of $u$, and call 
it the {\it constant term} of $u$.    

\item[$iii)$] The set of all the elements of $\rg$ 
with no constant term 
will be denoted by $V_{G, R}$, or simply by $V_G$ 
if the base ring $R$ is clear in the context.  

\item[$iv)$] When $R$ is an integral domain, 
by the {\it characteristic} of $R$ (denoted by $char.\,R$) we mean 
the {\it characteristic} of the field of fractions of $R$. 
\end{enumerate}

Next, we start with the following equivalent 
formulation of Problem \ref{MainProb} for  
the group algebras of finite groups. 

\begin{propo}\label{equprop1}
Let $R$ be any commutative ring and $G$ a finite group. 
Then $\vg$ is a Mathieu subspace of any fixed type of $R[G]$ 
iff all elements of $\sqrt{\vg}$ are nilpotent.  
\end{propo}
  
\pf First, it is easy to see that the $(\Leftarrow)$ part 
follows directly from the assumption and 
Definition \ref{Def-MS}. 

For the $(\Rightarrow)$ part, here we only give a proof 
for the left Mathieu subspace case.  
The proofs of the other three cases 
are similar. 
 
Assume that $\vg$ is a left Mathieu subspace  
and let $u\in \sqrt{\vg}$. Replacing $u$ by a positive 
power of $u$, if necessary, we may assume that 
$u^m\in \vg$ for all $m\ge 1$.

Now, since $G$ is finite, by 
Definition \ref{Def-MS} there exists $N\ge 1$ 
such that $g^{-1} u^m\in \vg$ for all $g\in G$ and $m\ge N$. 
In particular, for each $g\in G$, the constant term of 
$g^{-1}u^N$, which is the same as the coefficient of $g$ in $u^N$, 
is equal to $0$, whence $u^N=0$, i.e., $u$ is nilpotent.  

Another way to show the $(\Rightarrow)$ part is as follows. 

Assume otherwise and let $u\in \sqrt{\vg}$ such that $u^m\neq 0$ 
for all $m\ge 1$. Since $G$ is finite, there exists $g\in G$ 
such that the coefficient of $g$ in $u^m$ is nonzero for 
infinitely many $m\ge 1$. Then the constant term of 
$g^{-1}u^m$ is nonzero for infinitely many $m\ge 1$. 
Then by Definition \ref{Def-MS} $\vg$ is not 
a Mathieu subspace of $R[G]$, 
which is a contradiction.
\epfv

Two immediate consequences of Proposition \ref{equprop1} 
are the following two corollaries.

\begin{corol}\label{AllType}
Let $R$ and $G$ be as in Proposition \ref{equprop1}. 
Then $\vg$ is a Mathieu subspace of any fixed type of $R[G]$ 
iff $\vg$ is a $($two-sided$)$ Mathieu subspace of $R[G]$. 
\end{corol}

Therefore, throughout the rest of this paper 
we may and will focus only on the two-sided case.

\begin{corol}\label{OneWayIdem}
Let $R$ and $G$ be as in Proposition \ref{equprop1}. 
Assume that $\vg$ is a Mathieu subspace of $R[G]$. Then  
$\vg$ contains no nonzero idempotent 
of $\rg$. 
\end{corol}
\pf Assume otherwise. Let $e\in \vg$ be a 
nonzero idempotent, i.e., $e^2=e\ne 0$. 
Then for any $m\ge 1$, we have $e^m=e\in\vg$, 
whence $e\in\sqrt \vg$. But, since $e$ is 
clearly not nilpotent, by 
Proposition \ref{equprop1} $\vg$ is 
not a Mathieu subspace of $\rg$, 
which is a contradiction.
\epfv

When the base ring $R$ is a field, 
we show next that the converse of 
Corollary \ref{OneWayIdem} actually 
also holds. 

\begin{propo}\label{equprop2}
Let $K$ be a field and $G$ a finite group. 
Then $\vg$ is a Mathieu subspace of $\kg$ iff 
$\vg$ contains no nonzero idempotent of $\kg$.   
\end{propo} 

\pf The $(\Rightarrow)$ part 
is a special case of Corollary \ref{OneWayIdem}. 
To show the $(\Leftarrow)$ part, note that 
$\kg$ is algebraic over $K$, since 
it is of finite dimension over $K$.
In particular, the radical $\sqrt\vg$ of 
$\vg$ is algebraic over $K$.  
Then by Theorem \ref{CharByIdem}, 
$\vg$ is a Mathieu subspace of $\kg$. 
\epfv

Next, we show that Problem \ref{MainProb} 
can be solved for the group algebras 
of all finite groups $G$ over integral domains 
$R$ such that $char.\,R=0$ or $char.\,R=p>|G|$. 

\begin{theo}\label{GeneralCase}
Let $G$ be a  finite group  and $R$ an 
integral domain such that $char.\,R=0$ or $char.\,R=p>|G|$. 
Then $\vg$ is a Mathieu subspace of $\rg$. 
\end{theo}

\pf Let $u\in \sqrt\vg$. Then by Proposition \ref{equprop1} 
it suffices to show that $u$ is nilpotent.  
Note that by replacing $u$ by a positive power of $u$,  
if necessary, we may assume   
$u^m\in \vg$, i.e., $\tr(u^m)=0$, 
for all $m\ge 1$.  

Let $\mu:\rg \rightarrow {\rm End}_R(\rg)$ be the  
$R$-algebra homomorphism which maps 
each $v\in \rg$ to the $R$-endomorphism 
$m_v\in {\rm End}_R(\rg)$ defined by 
the left multiplication by $v$ on $\rg$. 
Then it is easy to check that 
for any $v\in \rg$, the trace of the 
linear map $\mu(v)=m_v$ is equal to $|G| \tr(v)$.  
Consequently, for the $u\in \sqrt\vg$ fixed at the 
beginning and any $m\ge 1$,  the trace of the 
$m$-th power $(\mu(u))^m=\mu(u^m)$ 
of the linear transformation $\mu(u)$ 
is equal to zero. 

On the other hand, since $char.\,R=0$ or $char.\,R=p>|G|$, 
it is well-known in linear algebra that in this case 
the linear transformation $\mu(u)$ must be nilpotent, 
i.e., $(\mu(u))^m=\mu(u^m)=0$ for $m\gg 0$.  
Since $\mu$ is clearly injective 
(e.g., by applying $\mu(v)$ to $1\in \rg$ 
for all $v\in \rg$), we also have $u^m=0$ 
when $m\gg 0$, i.e., $u$ is nilpotent, 
as desired.  
\epfv

One remark on Theorem \ref{GeneralCase} is that when 
the conditions $char.\,R=0$ and $char.\,R=p>|G|$ 
fail, i.e., when $0<char.\,R=p \le |G|$, 
the situation for Problem \ref{MainProb} 
becomes much more complicated.

For instance, as shown by the next lemma and also by  
Theorem \ref{MainAbThm} in Section \ref{S4},  
the magic condition $p \nmid |G|$ for  
the theory of group algebras $\rg$ 
of finite groups $G$ 
(e.g., see \cite{P}) does not resolve 
the difficulty completely for 
Problem \ref{MainProb}. 

\begin{lemma}\label{UniCounter} 
Let $G$ be any finite group with $|G|\ge 2$, and 
$R$ an integral domain of $char.\,R=p>0$. 
Assume $p\,|\,(|G|-1)$ $($hence, $p\nmid |G|$$)$. 
Then $\vg$ is not a Mathieu subspace of $\rg$.  
\end{lemma} 
\pf Let $u=-\sum_{g\in G\backslash \{1_G\}}g\in \vg$ 
and $v=1_G-u=1-u$. Note that $v$ is the sum of all the 
distinct elements of $G$ in $\rg$.  
Hence, for any $g\in G$, we have 
$vg=gv=v$. Consequently, we have 
$v^2=|G|v$, which in terms of $u$ is the 
same as 
\begin{align*}
(1-u)^2=1-2u+u^2=|G|(1-u).
\end{align*}
Solving $u^2$ from the equation above, we get
\begin{align} 
u^2&=(|G|-1)-(|G|-2)u. \label{UniCounter-pe1}
\end{align} 

Since $p\,|\,(|G|-1)$, we have $(|G|-1)=0$ and 
$(|G|-2)=-1$. Then by Eq.\,(\ref{UniCounter-pe1}), 
we have $u^2=u$. Since $u\ne 0$, 
by Corollary \ref{OneWayIdem} $\vg$ is not 
a Mathieu subspace of $\rg$. 
\epfv

Next, we show the following lemma  
that will be needed later. 

\begin{lemma} \label{GrpExt1} 
Let $R$ be any commutative ring and $G$ 
any group $($not necessarily finite$)$. 
Assume that $\vg$ is a Mathieu subspace of $\rg$. 
Then for each subgroup $H$ of $G$,  
$V_H$ is a Mathieu subspace of $R[H]$. 
\end{lemma} 
\pf Assume otherwise. Let $H$ be a subgroup of 
$G$ such that $V_H$ is not a Mathieu subspace of 
$R[H]$. Then by Definition \ref{Def-MS} and 
the definition of $V_H$, 
there exist $u, v\in R[H]$ such that 
$\tr(u^m)=0$ for all $m\ge 1$, but 
$\tr(u^mv)\ne 0$ for infinitely 
many $m\ge 1$. 

Since $R[H]\subseteq \rg$, we have $u, v\in \rg$, 
and $u^m \in \vg$ for all $m\ge 1$, but 
$u^m v \not \in \vg$ for infinitely 
many $m\ge 1$. Hence, $\vg$ is not 
a Mathieu subspace of $\rg$,  
which is a contradiction. 
\epfv

\begin{corol} \label{GrpExt} 
Let $R$ and $G$ be as in Lemma \ref{GrpExt1} 
and $H$ a subgroup of $G$. 
Assume that $V_H$ is not a 
Mathieu subspace of $R[H]$.  
Then $\vg$ is not a Mathieu subspace 
of $\rg$. 
\end{corol} 

As an application of Lemma \ref{GrpExt1} or 
Corollary \ref{GrpExt}, we derive the following 
necessary condition for $\vg$ to be a Mathieu subspace of 
$\rg$ over integral domains $R$ of positive characteristic. 

\begin{propo}\label{q-Sylow}
Let $R$ be an integral domain of characteristic $p>0$ 
and $G$ an arbitrary finite group. Write $|G|=p^r d$ for some 
$r\ge 0$ and $d\ge 1$ with $p\nmid d$.  
Assume that $R$ contains a primitive $d$-th root 
of unity and $\vg$ is a Mathieu subspace of $\rg$. 
Then for each prime divisor $q$ of $|G|$, 
we have $p \ge q$. 
\end{propo} 

\pf Assume otherwise and let $q$ be a prime divisor 
of $|G|$ such that $p<q$. Then we have $q\,|\, d$, 
whence $R$ also contains a primitive $q$-th root 
of unity.   

Write $|G|=q^s n$ with $s, n\ge 1$ such that 
$q \nmid n$. Then by the well-known 
Sylow's theorem in the theory of finite groups 
(e.g., see p.\,$105$, Theorem $2.11.7$ 
in \cite{He}), $G$ has at least one 
$q$-Sylow subgroup $H$, i.e., 
a subgroup $H$ of $G$ with $|H|=q^s$. 

Now, pick up any non-identity element $h\in H$. 
Then $h$ has order $q^k$ for some $1\le k\le r$. 
Let $g=h$ if $k=1$; and $g=h^{k-1}$  
if $k\ge 2$. Then $g$ has order $q$ 
and hence, generates a cyclic subgroup 
$C_q$ of $G$ of order $|C_q|=q$.  
Then by Theorem \ref{MainAbThm} 
to be proved in Section \ref{S4}, 
$V_{C_q}$ is not a Mathieu subspace of $R[C_q]$. 
Hence, by Corollary \ref{GrpExt} 
$\vg$ is not a Mathieu subspace of $\rg$ either, 
which is a contradiction. 
\epfv

Finally, we point out that when the finite 
group $G$ in Proposition \ref{q-Sylow} 
is abelian, a much stronger condition will 
be given in Theorem \ref{MainAbThm}  
of the next section.

\renewcommand{\theequation}{\thesection.\arabic{equation}}
\renewcommand{\therema}{\thesection.\arabic{rema}}
\setcounter{equation}{0}
\setcounter{rema}{0}
 
\section{\bf The Case for Finite Abelian Groups}\label{S4}

In this section, we study Problem \ref{MainProb} for 
finite abelian groups over certain integral domains. 
The main result of this section is 
the following theorem. 

\begin{theo}\label{MainAbThm}
Let $R$ be an integral domain of characteristic $p>0$, 
and $G$ a finite abelian group with $|G|=p^rd$ for some 
$r\ge 0$ and $d\ge 1$ with $p\nmid d$. 
Assume that $R$ contains a primitive $d$-th root of unity. 
Then $\vg$ is a Mathieu subspace of $\rg$ iff $p>d$. 
\end{theo}
 
Two remarks on Theorem \ref{MainAbThm} are as follows. 
\vspace{2mm}

First, when the integral domain $R$ has $char.\,R=0$ 
(or $char.\,R=p>|G|$), Problem \ref{MainProb} has been 
solved by Theorem \ref{GeneralCase}, 
together with which Theorem \ref{MainAbThm} 
provides a complete solution  
of Problem \ref{MainProb} for the group algebras of 
all finite abelian groups when the base integral 
domain $R$ satisfies the primitive root of unity condition in 
Theorem \ref{MainAbThm}, e.g., when 
$R$ is an algebraically closed field. 

Second, from the example below we see that the $d$-th primitive 
root of unity condition on the integral domain $R$ in 
Theorem \ref{MainAbThm} is necessary. 

\begin{exam}\label{F3Z5} 
Let $\mathbb F_3$ be the field with three elements. 
Note that $\mathbb F_3$ obviously does not 
contain any primitive $5$th root of unity. 
But, $V_{\bZ_5}$ is a Mathieu subspace of 
$\hspace{.2mm}\mathbb F_3[\bZ_5]$, although 
$char.\,\mathbb F_3=3<d=5$.
\end{exam}
\pf Assume otherwise. Then by Proposition \ref{equprop2},  
there exists a nonzero idempotent 
$f\in V_{\bZ_{5}}$. 
By identifying the group algebra 
$\mathbb F_3[\bZ_5]$ with the quotient algebra 
$\mathbb F_3[t]/(t^5-1)$ of the polynomial 
algebra $\mathbb F_3[t]$ in one variable $t$, 
we may write $f=c_1t+c_2t^2+c_3t^3+c_4t^4$.
Then it is easy to check that the following 
equations hold: 
\begin{align*}
\tr(f^2)=2(c_1c_4+c_2c_3), \\ 
f^3=c_1^3 t^3 +c_2^3 t+c_3^3 t^4+ c_4^3t^2.
\end{align*}

Since $f^2=f^3=f\in V_{\bZ_5}$, hence we also have  
\begin{align}
c_1c_4 =-c_2c_3, \label{F3Z5-pe3}
\end{align}
\begin{align}
c_1=c_2^3;\,\, \,\, c_2=c_4^3;\,\, \,\, c_3
=c_1^3;\,\, \,\, c_4=c_3^3. \label{F3Z5-pe4}
\end{align}

From the four equations in Eq.\,(\ref{F3Z5-pe4}), 
it is easy to see that if one of the $c_i$'s is equal to zero, 
then so are all the $c_i$'s. Since $f\ne 0$, we see that 
all the $c_i$'s are nonzero. 

By combining equations in Eqs.\,(\ref{F3Z5-pe3})-(\ref{F3Z5-pe4}),  
it is also easy to see that $(c_2c_3)^3 =-(c_2c_3)$,  
whence $(c_2c_3)^2 =-1$.  However, the base field 
$\mathbb F_3$ contains no square root of $-1$.  
Hence, we get a contradiction. 
\epfv

Next, we will devote the rest of this section to give 
a proof for Theorem \ref{MainAbThm}. First, 
we need to show the following reduction lemma.  

\begin{lemma}\label{p-cross} 
Let $R$ be an integral domain of characteristic $p>0$ 
and $H$ a finite abelian group. Let $q=p^r$ for some $r\ge 1$ 
and $G=H\times \bZ_{q}$.  
Then $V_H$ is a Mathieu subspace of $R[H]$ 
iff $\vg$ is a Mathieu subspace of $\rg$. 
\end{lemma}

\pf For convenience, we identify $\bZ_q$ 
with the multiplicative cyclic group $C_q$  
with $q$-element. We also identify 
$H$ and $C_q$ with the subgroups $H\times\{1_{C_q}\}$ 
and $\{1_H\}\times C_{q}$ of $G$, 
respectively. 

Under these identifications, 
$G$ is also the inner product of its subgroups 
$H$ and $C_q$, and the group algebras 
$R[H]$ and $R[C_q]$ become subalgebras 
of $\rg$. Then the $(\Leftarrow)$ part
of the lemma follows immediately 
from Lemma \ref{GrpExt1}. 

To show the $(\Rightarrow)$ part, pick up any 
$u\in \sqrt{\vg}$. Then by Proposition \ref{equprop1},  
it suffices to show that $u$ is nilpotent.  To do so, 
replacing $u$ by a positive power of $u$, if necessary, 
we assume that $u^m\in \vg$ for all $m\ge 1$.  

Write $u=\sum_{s\in C_q} \alpha_s s$ with 
$\alpha_s\in R[H]$ for each $s\in C_q$. 
Note that for any $k\ge 1$ and $s\in C_q$, 
we have $s^{q^k}=1_{C_q}$, since $|C_q|=q$.    
Then by the conditions that 
$char.\,R=p>0$ and $q$ is a positive power of $p$, 
for any $k\ge 1$ we also have   
\begin{align*}
u^{q^k} = \sum_{s\in C_q} \alpha_s^{q^k} s^{q^k} 
= \sum_{s\in C_q} \alpha_s^{q^k} \in R[H].  
\end{align*} 

Moreover, since $u^m\in \vg$ for all $m\ge 1$, 
we have $(u^q)^k=u^{q^k}\in R[H] \cap \vg=V_H$ 
for all $k\ge 1$, whence $u^q \in \sqrt{V_H}$. 
Since by assumption $V_H$ is a Mathieu subspace 
of $R[H]$, applying Proposition \ref{equprop1} 
to the group algebra $R[H]$ we see that $u^q$ 
is nilpotent, whence so is $u$.  
\epfv

Next, let's recall the following well-known 
fundamental theorem of finite abelian groups.

\begin{theo}\label{G-cross}
Any finite abelian group can be written as a direct product of 
cyclic groups whose orders are powers of primes. 
\end{theo} 

For the proof of the theorem above, see any abstract algebra 
text book (e.g., see Th.2.2, Ch.II, \cite{H}). 

Note that by applying Theorem \ref{G-cross} and 
Lemma \ref{p-cross} (inductively), 
it is easy to see that we may actually 
assume that the exponent $r$ in 
Theorem \ref{MainAbThm} is equal to zero, 
i.e., it suffices to show the 
following lemma.
 
\begin{lemma}\label{SpecialCase} 
Let $G$ be a finite abelian group and 
$R$ an integral domain of characteristic 
$p>0$ such that $p\nmid d\!:=|G|$.
Assume that $R$ contains a primitive 
$d$-th root of unity. 
Then $\vg$ is a Mathieu 
subspace of $\rg$ iff $p>d=|G|$.  
\end{lemma}

From now on and throughout the rest of 
this section, we let $G$ and $R$ be as 
in the lemma above. 

Note first that when $d=|G|=1$, 
we have $\vg=\{0\}$, which is obviously 
a Mathieu subspace 
of $\rg$. Hence, Lemma \ref{SpecialCase} holds 
in this trivial case. So we will  
assume $d=|G|\ge 2$. 

Note also that by Theorem \ref{G-cross}, 
we may (and will) further assume that 
the abelian group $G$ is given by  
\begin{align} 
G=\bZ_{d_1}\times \bZ_{d_2}\times 
\cdots \times \bZ_{d_n} \label{GenericAbG}  
\end{align}
for some $n\ge 1$ and $d_i\ge 2$ $(1\le i\le n)$. 

But, here we do not need to assume that 
the integers $d_i\ge 2$ $(1\le i\le n)$ 
are powers of primes. 

In order to study the group algebra $\rg$ of $G$ 
in Eq.\,(\ref{GenericAbG}), we need to write 
the factor groups $\bZ_{d_i}$ 
$(\lin)$ in Eq.\,(\ref{GenericAbG}) 
as multiplicative groups $H_i$ 
with a fixed generator $e_i\in H_i$, 
i.e., for each $\lin$, we let 
\begin{align}
H_i=\{ e_i^k \,|\, 0\le k\le d_i-1\} \simeq \bZ_{d_i}.   
\end{align}

For convenience, for each $\lin$,  
we also identify $H_i$ (implicitly) with the 
subgroup of $G$ in Eq.\,(\ref{GenericAbG}) 
consisting of all the $n$-tuples whose 
$j$-th $(j\ne i)$ component  
being the identity element of 
$H_j\simeq\bZ_{d_j}$. 
Note that under this identification, 
we have $H_i\subset G$, whence $G$ is also 
the {\it inner product} of 
the subgroups $H_i$ $(\lin)$, 
i.e., with the abusive notations fixed above, 
we have
\begin{align}
G=H_1\cdot H_2 \cdots H_n = H_1 \times H_2 \times 
\cdots \times H_n  \label{GenericAbG-2}  
\end{align}

Furthermore, we also need to introduce 
the following two sets:
 \begin{align}
D\!:=&\{\beta=(\beta_1, \beta_2, ..., \beta_n) \in \bN^n \,|\, 
0 \le \beta_i\le d_i-1  \} \label{Def-D} \\
S\!:=&\left\{a=(a_1, a_2, ... , a_n)\in R^n \,|\,\, 
a_i^{d_i}=1  \right\}.  \label{Def-S}  
\end{align} 

Note that since $R$ contains a primitive $d$-th root of unity,  
$R$ also contains a primitive $d_i$-{\rm th} $(\lin)$ 
root of unity, since $d_i\,|\, d$. Then from Eqs.\,(\ref{Def-D}) 
and (\ref{Def-S}), we have $|S|=d=|D|=|G|$. 
\vspace{2mm} 

Next, with the notations fixed above we give an 
equivalent formulation of Lemma \ref{SpecialCase} 
in terms of the polynomial algebra $\rz$ over $R$ 
in $n$ variables $z\!:=(z_1, z_2, ... , z_n)$. 

First, we define and consider the following $R$-linear functional: 
\begin{eqnarray} 
\cL\!: \rz\, & \to & \quad R \label{Def-cL}  \\
f  \quad &  \to & \sum_{a\in S} f(a). \nno
\end{eqnarray}  

\begin{lemma}\label{SumLemma} 
Let $G$ and $R$ be fixed as above. 
Then for any $\alpha \in D$, we have  
\begin{align}\label{SumLemma-e1} 
\cL (z^\alpha)=\begin{cases} 
d   &\mbox{ if } \alpha=0;\\
0   &\mbox{ if } \alpha \neq 0. 
\end{cases} 
\end{align}  
\end{lemma}
\pf If $\alpha=0$, then $\cL (z^\alpha)=\sum_{a\in S}1=|S|=d$. 
So we let $\alpha\ne 0$. Without losing any generality, 
we assume that the first component of $\alpha$ is nonzero, 
and denote it by $k$ (for short). 

Let $\xi_1$ be a primitive $d_1$-th 
root of unity in $R$. Then we have $\xi_1^k \ne 1$, 
since $1\le k \le d_1-1$.  
Note that for each root $1 \ne r \in R$ of the polynomial 
$z_1^{d_1}-1\in R[z_1]$, $r$ is also a root of the polynomial 
$\sum_{\ell=0}^{d_1-1} z_1^\ell$, for  
$z_1^{d_1}-1=(z_1-1)\sum_{\ell=0}^{d_1-1} z_1^\ell$. 
Therefore, for the fixed primitive $d_1$-th 
root of unity $\xi_1\in R$, we have 
\begin{align}
\sum_{\ell=0}^{d_1-1} (\xi_1^\ell)^k=
\sum_{\ell=0}^{d_1-1} (\xi_1^k)^\ell=0.  
\label{SumLemma-1pe}
\end{align}

Now, for each $\lin$, set 
$C_i\!:=\{\xi_i^\ell \,|\, 0 \le \ell\le d_i-1\}$, 
where $\xi_i$ is any fixed primitive $d_i$-th root 
of unity in $R$. Then from the definition of 
the set $S$ in Eq.\,(\ref{Def-S}), we have 
$S=C_1\times C_2 \times \cdots \times C_n$.  
By taking the sum 
$\cL(z^\alpha)=\sum_{a\in S} a^\alpha$ 
first over the set $C_1$, 
it follows immediately from  
Eq.\,(\ref{SumLemma-1pe}) 
that $\cL(z^\alpha)=0$.  
\epfv

Next, we define the following $R$-algebra 
homomorphism: 
\begin{align}
\varphi: \rz & \to \rg \label{Def-varphi} \\ 
 z_i  \,\, & \to \,\,\, e_i.  \nno
\end{align} 

Note that the kernel of the $R$-algebra 
homomorphism $\varphi$ above is the ideal 
of $\rz$ generated by the polynomials 
$z_i^{d_i}-1$ $(\lin)$. We will denote 
this ideal by $I_{\vec d}$, where 
$\vec{d}$  stands for the $n$-tuple 
$(d_1, d_2, ..., d_n)$.  
 
The pre-image of $\vg \subset \rg$ 
under the linear map $\varphi$  
is given by the following lemma.  

\begin{lemma}\label{CrucialLemma} 
With the setting above,  we have  
\begin{align}\label{CrucialLemma-e1} 
\varphi^{-1}(\vg)=\Ker \, \cL.  
\end{align}  
\end{lemma}

\pf First, let $V_0$ be the $R$-subspace of $\rz$ spanned 
by $z^\alpha$ $(0\ne \alpha \in D)$ and 
$V\!:=R\cdot 1\oplus V_0$. Then by the definition of $\varphi$ 
in Eq.\,(\ref{Def-varphi}), it is easy to see that we have 
\begin{align}\label{CrucialLemma-pe1}
\varphi^{-1}(\vg)=\big\{ f \in \rz \,\big|\, 
f \equiv r \,\,(\mbox{\rm mod\,} I_{\vec d}\hspace{.02mm}) 
\,\, \mbox{ for some } r\in V_0 \big\}.  
\end{align}

Therefore, it suffices to show that $\Ker \, \cL$ 
coincides with the set on the right-hand side 
of the equation above.

Now, let $f\in \rz$. Then there exists a unique $r\in V$ 
such that $f\equiv r \,\,(\mbox{\rm mod\,} I_{\vec d}\hspace{.02mm})$. 
By Eq.\,(\ref{CrucialLemma-pe1}) we have  
\begin{align}
f\in \varphi^{-1}(\vg) \,\,\, \Leftrightarrow  \,\,\, r\in V_0. 
\label{CrucialLemma-pe3}  
\end{align} 

Furthermore, since $S$ is the zero-set 
of the ideal $I_{\vec d}$ in $R^n$, 
we have $f(a)=r(a)$ for all $a\in S$.
In particular, we have 
$\cL(f)=\cL(r)$ and hence,  
\begin{align}
f\in \Ker \cL \,\,\, \Leftrightarrow  \,\,\, r\in \Ker \cL. \label{CrucialLemma-pe2}  
\end{align}

Write $r(z)=\sum_{\alpha\in D} c_\alpha z^\alpha$. 
Then by Eq.\,(\ref{SumLemma-e1}) we have 
\begin{align*}
\cL(r)=\cL(c_0)+\sum_{0\ne \alpha\in D} c_\alpha \cL(z^\alpha)
=dc_0.
\end{align*}

Since $p\nmid d$, we see that $r\in \Ker \cL$ 
iff $c_0=0$ iff $r\in V_0$. Then by the equivalences in 
Eqs.\,(\ref{CrucialLemma-pe3}) and (\ref{CrucialLemma-pe2}), 
we have that $f\in \varphi^{-1}(\vg)$ iff $f\in \Ker \cL$, 
whence the lemma follows. 
\epfv

Finally, we can give a proof for Lemma \ref{SpecialCase} as follows,  
from which the proof of the main result Theorem \ref{MainAbThm}  
will be completed.  \\

\underline{\it Proof of Lemma \ref{SpecialCase}:} 
Note that the $(\Leftarrow)$ part of the lemma follows 
directly from Theorem \ref{GeneralCase}, 
which actually does not need the primitive 
root of unity condition on $R$ in the lemma. 
But, with the primitive 
root of unity condition on $R$ it also follows 
from the arguments below. 

First, we consider the $R$-homomorphism 
$\varphi: \rz\to \rg$ defined in 
Eq.\,(\ref{Def-varphi}). Note that $\varphi$ 
is surjective with the kernel $I_{\vec d}$. 
Hence, from Eq.\,(\ref{CrucialLemma-e1})  
we have $I_{\vec d} \subseteq \Ker\cL$ and 
$\varphi(\Ker\cL)=\vg$. 

Therefore, we may identify $\rg$ with the quotient 
algebra $\rz/I_{\vec d}$, and $V_G$ with 
$\Ker\cL /I_{\vec d}$. Via these identifications and  
by Proposition \ref{quotient-I}, we have that 
$\vg$ is a Mathieu subspace of $\rg$, iff   
$\Ker\cL$ is a Mathieu subspace of 
the polynomial algebra $\rz$. 

Second, by applying Proposition \ref{FiniteCase2} 
to the set $S$ in Eq.\,(\ref{Def-S}) with $c_i=1$ 
$(1\le i\le d)$, we have that $\Ker\cL$ is 
a Mathieu subspace of $\rz$, iff for  
any non-empty subset $J\subseteq\{1, 2, ..., d\}$, 
the cardinal number $|J|\neq 0$ in $R$, 
i.e., $|J|\not\equiv 0\mod p$.  
Furthermore, it is easy to see that the latter 
property holds iff $p>d=|G|$.  

Finally, by combining the three equivalences above, 
we see that the lemma follows. 
\epfv

\renewcommand{\theequation}{\thesection.\arabic{equation}}
\renewcommand{\therema}{\thesection.\arabic{rema}}
\setcounter{equation}{0}
\setcounter{rema}{0}
 
\section{\bf The Case for the Group Algebra $R[\bZ^n]$ 
with $char.\,R=p>0$}\label{S5}

In this section, we show that Problem \ref{MainProb} 
has a negative answer for the group algebras of 
the free abelian groups $\bZ^n$ $(n\ge 1)$ 
over all integral domains $R$ of 
positive characteristics. More precisely, 
we have the following proposition. 

\begin{propo} \label{Not-MS} 
For any integral domain $R$ of $char.\,R=p>0$, 
$V_{\bZ^n}$ is not a Mathieu subspace of the 
group algebra $R[\bZ^n]$. 
\end{propo}

Note that under the natural identification 
$R[\bZ^n]\simeq R[z^{-1},z]$ 
(the Laurent polynomial algebra in $n$ variables 
$z=(z_1, z_2, ..., z_n)$ over $R$),   
the proposition above is equivalent to saying 
that for any integral domain $R$ of $char.\,R=p>0$, 
the subspace $V$ of all the Laurent polynomials in 
$R[z^{-1}, z]$ with no constant term does not form 
a Mathieu subspace of the Laurent polynomial algebra 
$R[z^{-1}, z]$. In particular, it follows that 
the Duistermaat-van der Kallen Theorem, 
Theorem \ref{ThmDK}, cannot be generalized 
to any field of characteristic $p>0$.    

To show Proposition \ref{Not-MS}, 
note first that we may identify  
$\bZ$ as the subgroup of $\bZ^n$
consisting of all the elements 
$(a, a, ..., a)\in \bZ^n$ 
with $a\in \bZ$. Then by 
Corollary \ref{GrpExt}, 
we may actually assume $n=1$. 
Furthermore, via the identification 
$R[\bZ]\simeq R[z, z^{-1}]$ mentioned 
above, it will be enough to show the 
following lemma. The example in the lemma 
was suggested to the authors 
by Arno van den Essen. 

\begin{lemma} \label{ExampleArno} 
Let $p$ be a prime and $z$ a free variable.  
Set $f\!:=z^{-1}+z^{p-1} \in \bZ_p[z^{-1}, z]$. 
Then the following two statements hold:  
\begin{enumerate}
\item[$i)$]  $\tr(f^m)=0$ for all $m\ge 1$; 
\item[$ii)$] $\tr\big(z^{-1}f^{p^k-1}\big)=(-1)^{p^{k-1}}$ 
            for all $k\ge 1$. 
\end{enumerate} 
\end{lemma}

In order to prove the lemma above, we need first to 
show the following lemma. 

\begin{lemma} \label{p-lemma} 
For any prime number $p>0$, the following statements hold.  
\begin{enumerate}
\item[$i)$] For any $k, a \in \bN$ such that $k\ge 1$ and $a\le p^k-1$, 
   we have 
\begin{align}\label{p-lemma-e2}
\binom{p^k-1}{a}\equiv (-1)^a\mod p.
\end{align} 
\item[$ii)$] For any integer $b\ge 1$, we have 
\begin{align}\label{p-lemma-e1}
\binom{bp}{b}\equiv 0 \mod p.
\end{align} 
\end{enumerate}
\end{lemma}
\pf $i)$ Let $x$ be a free variable. 
We consider the polynomial $(x-1)^{p^k-1}$ 
in the rational function field $\bZ_p(x)$,  
for which we have the following two equations: 
\begin{align}
(1-x)^{p^k-1}&=\sum_{a=0}^{p^k-1} (-1)^a  
\binom{p^k-1}a x^a,  \label{p-lemma-pe2} \\
(1-x)^{p^k-1}&=\frac{(1-x)^{p^k}}{1-x}
=\frac{1-x^{p^k}}{1-x}=\sum_{a=0}^{p^k-1} x^a. \label{p-lemma-pe3}
\end{align} 
Note that Eq.\,(\ref{p-lemma-pe3}) above also holds 
for the case $p=2$, since $1=-1$ in $\bZ_2$. 
Now, by comparing the coefficients of $x^a$ 
in the polynomials on the right-hand sides of 
Eqs.\,(\ref{p-lemma-pe2}) and (\ref{p-lemma-pe3}), 
we see that $i)$ follows.  

$ii)$ Write $b=p^rn$ for some $r\ge 0$ and $n\ge 1$ 
such that $p \nmid n$. In particular, we have  
$p^{r+1} \nmid\, b$. 

We consider the polynomial $(x+1)^{bp}\in \bZ_p[x]$. 
Note that the coefficient of $x^b$ in $(x+1)^{bp}$ 
is equal to $\binom{bp}{b}$. On the other hand, 
we also have 
\begin{align*}
(x+1)^{bp}&=(x+1)^{np^{r+1}}=(x^{p^{r+1}}+1)^n.
\end{align*} 

Now, assume that $\binom{bp}{b}\not\equiv 0\mod p$. 
Then by the equation above, 
$x^b$ appears in the polynomial 
$(x^{p^{r+1}}+1)^n$ with a nonzero 
coefficient, whence $b=p^{r+1}k$ 
for some $1\le k\le n$. But this implies 
$p^{r+1}\,|\, b$, which is a contradiction.     
\epfv

\underline{\it Proof of Lemma \ref{ExampleArno}:} \, 
$i)$ Since $f=z^{-1}+z^{p-1}$, the constant term 
of $f^m$ $(m\ge 1)$ is given by the sum of $\binom mb$ 
for  all the integers $0\le b\le m$ such that 
$-(m-b)+b(p-1)=0$, which is the same as $m=bp$. 
Therefore, there is at most one such an integer $b$, 
which is $m/p$ if (and only if) $p\,|\, m$. Hence   
we have    
\begin{align}
\tr(f^m)=\begin{cases}
\binom{bp}{b} &\mbox{ if } p\, | \, m \mbox{ and } b=m/p; \\
\,\,\,\, 0 &\mbox{ if } p \nmid m. 
\end{cases}
\end{align} 
Then from the equation above and Eq.\,(\ref{p-lemma-e1}), 
we see that $i)$ follows.   

$ii)$ By a similar argument as in $i)$, it is easy to  
check that for any $k\ge 1$, the coefficient of 
$z$ in $f^{p^k-1}$ is given by $\binom{p^k-1}{p^{k-1}}$, 
which by Eq.\,(\ref{p-lemma-e2}) is equal to $(-1)^{p^{k-1}}$. 
Hence, we have $\tr(z^{-1}f^{p^k-1})=(-1)^{p^{k-1}}$  
for all $k\ge 1$, i.e., $ii)$ holds. 
\epfv

{\bf Acknowledgments}\,\,The authors are very grateful to  
Professor Arno van den Essen for suggesting the example 
in Lemma \ref{ExampleArno}.

\end{document}